\documentclass[11pt]{amsart}
\usepackage{amsmath,amssymb}
\newtheorem{theorem}{Theorem}
\newtheorem{proposition}[theorem]{Proposition}

\newtheorem{lemma}[theorem]{Lemma}

\title[Flows by powers of the Gaussian curvature]{Asymptotic behavior of flows by powers of the Gaussian curvature}

\author{Simon Brendle}
\address{Department of Mathematics, Columbia University, 2990 Broadway, New York, NY 10027, USA.}
\email{simon.brendle@columbia.edu}
\author{Kyeongsu Choi}
\address{Department of Mathematics, Columbia University, 2990 Broadway, New York, NY 10027, USA.}
\email{kschoi@math.columbia.edu}
\author{Panagiota Daskalopoulos}
\address{Department of Mathematics, Columbia University, 2990 Broadway, New York, NY 10027, USA.}
\email{pdaskalo@math.columbia.edu}
\thanks{The first author was supported in part by the National Science Foundation under grant DMS-1649174. The third author was supported in part by the National Science Foundation under grant DMS-1600658.}

\begin{document}

\begin{abstract}
We consider a one-parameter family of strictly convex hypersurfaces in $\mathbb{R}^{n+1}$ moving with speed $- K^\alpha \nu$, where $\nu$ denotes the outward-pointing unit normal vector and $\alpha \geq \frac{1}{n+2}$. For $\alpha > \frac{1}{n+2}$, we show that the flow converges to a round sphere after rescaling. In the affine invariant case $\alpha=\frac{1}{n+2}$, our arguments give an alternative proof of the fact that the flow converges to an ellipsoid after rescaling.
\end{abstract}

\maketitle

\section{Introduction}

Parabolic flows for hypersurfaces play an important role in differential geometry. One fundamental example is the flow by mean curvature (see \cite{Huisken}). In this paper, we consider flows where the speed is given by some power of the Gaussian curvature. More precisely, given an integer $n \geq 2$ and a real number $\alpha >0$, a one-parameter family of immersions $F:M^n \times [0,T) \to \mathbb{R}^{n+1}$ is a solution of the $\alpha$-Gauss curvature flow, if for 
each $t \in [0,T)$, $F(M^n,t)=\Sigma_t$ is a complete convex hypersurface in $\mathbb{R}^{n+1}$, and $F(\cdot,t)$ satisfies
\begin{equation*}
\frac{\partial}{\partial t}  F(p,t)=- K^\alpha(p,t) \nu(p,t). 
\end{equation*} 
Here, $K(p,t)$ and $\nu(p,t)$ are the Gauss curvature and the outward pointing unit normal vector of $\Sigma_t$ at the point $F(p,t)$, respectively. 

\begin{theorem}
\label{main.thm}
Let $\Sigma_t$ be a family of closed, strictly convex hypersurfaces in $\mathbb{R}^{n+1}$ moving with speed $-K^\alpha \nu$, where $\alpha \geq \frac{1}{n+2}$. Then either the hypersurfaces $\Sigma_t$ converge to a round sphere after rescaling, or we have $\alpha = \frac{1}{n+2}$ and the hypersurfaces $\Sigma_t$ converge to an ellipsoid after rescaling.
\end{theorem}

Flows by powers of the Gaussian curvature have been studied by many authors, starting with the seminal paper of Firey \cite{Firey} in 1974. Tso \cite{Tso} showed that the flow exists up to some maximal time, when the enclosed volume converges to $0$. In the special case $\alpha=\frac{1}{n}$, Chow \cite{Chow1} proved convergence to a round sphere. Moreover, Chow \cite{Chow2} obtained interesting Harnack inequalitites for flows by powers of the Gaussian curvature (see also \cite{Hamilton}). In the affine invariant case $\alpha=\frac{1}{n+2}$, Andrews \cite{Andrews1} showed that the flow converges to an ellipsoid. This result can alternatively be derived from a theorem of Calabi \cite{Calabi}, which asserts that the only self-similar solutions for $\alpha=\frac{1}{n+2}$ are ellipsoids. The arguments in \cite{Calabi} and \cite{Andrews1} rely crucially on the affine invariance of the equation and do not generalize to other exponents. In the special case of surfaces in $\mathbb{R}^3$ ($n=2$), Andrews \cite{Andrews2} proved that flow converges to a sphere when $\alpha=1$; this was later extended in \cite{Andrews-Chen} to the case $n=2$ and $\alpha \in [\frac{1}{2},2]$. The results in \cite{Andrews2} and \cite{Andrews-Chen} rely on an application of the maximum principle to a suitably chosen function of the curvature eigenvalues; these techniques do not appear to work in higher dimensions. However, it is known that the flow converges to a self-similar solution for every $n \geq 2$ and every $\alpha \geq \frac{1}{n+2}$. This was proved by Andrews \cite{Andrews3} for $\alpha \in [\frac{1}{n+2},\frac{1}{n}]$; by Guan and Ni \cite{Guan-Ni} for $\alpha=1$; and by Andrews, Guan, and Ni \cite{Andrews-Guan-Ni} for all $\alpha \in (\frac{1}{n+2},\infty)$. One of the key ingredients in these results is a monotonicity formula for an entropy functional. This monotonicity was discovered by Firey \cite{Firey} in the special case $\alpha=1$.

Thus, the problem can be reduced to the classification of self-similar solutions. In the affine invariant case $\alpha=\frac{1}{n+2}$, the self-similar solutions were already classified by Calabi \cite{Calabi}. In the special case when $\alpha \geq 1$ and the hypersurfaces are invariant under antipodal reflection, it was shown in \cite{Andrews-Guan-Ni} that the only self-similar solutions are round spheres. Very recently,  the case $\frac{1}{n} \leq \alpha < 1+\frac{1}{n}$ was solved in \cite{Choi} as part  of K. Choi's Ph.D. thesis. In particular, this includes the case $\alpha=1$ conjectured by Firey \cite{Firey}. 

Finally, we note that there is a substantial literature on other fully nonlinear parabolic flows for hypersurfaces (see e.g. \cite{Chow-Tsai}, \cite{Huisken-Polden}, \cite{Schnurer}) and for Riemannian metrics (cf. \cite{Chow-Hamilton}). In particular, Gerhardt \cite{Gerhardt} studied convex hypersurfaces moving outward with speed $K^\alpha \nu$, where $\alpha<0$. Note that, for $\alpha<0$, one can show using the method of moving planes that any convex hypersurface satisfying $K^\alpha = \langle x,\nu \rangle$ where $\alpha<0$ is a round sphere. Using an a priori estimate in \cite{Chow-Gulliver}, Gerhardt \cite{Gerhardt} proved that the flow converges to a round sphere after rescaling.

We now give an outline of the proof of Theorem \ref{main.thm}. In view of the discussion above, it suffices to classify all closed self-similar solutions to the flow. The self-similar solutions $\Sigma=F(M^n)$ satisfy the equation
\begin{align*}\label{eq:INT Shrinker}
K^{\alpha}=\langle F,\nu \rangle. \tag{$*_\alpha$}
\end{align*}
To classify the solutions of \eqref{eq:INT Shrinker}, we distinguish two cases. 

First, suppose that $\alpha \in [\frac{1}{n+2},\frac{1}{2}]$. In this case, we consider the quantity 
\[Z = K^\alpha \, \text{\rm tr}(b) - \frac{n\alpha-1}{2\alpha} |F|^2,\] 
where $b$ denotes the inverse of the second fundamental form. The motivation for the quantity $Z$ is that $Z$ is constant when $\alpha=\frac{1}{n+2}$ and $\Sigma$ is an ellipsoid. Indeed, if $\alpha=\frac{1}{n+2}$ and $\Sigma = \{x \in \mathbb{R}^{n+1}: \langle Sx,x \rangle = 1\}$ for some positive definite matrix $S$ with determinant $1$, then $K^{\frac{1}{n+2}} = \langle F,\nu \rangle$ and 
\[Z = K^{\frac{1}{n+2}} \, \text{\rm tr}(b) + |F|^2 = \text{\rm tr}(S^{-1}).\] 
Hence, in this case $Z$ is constant, and equals the sum of the squares of the semi-axes of the ellipsoid.

Suppose now that $\Sigma = F(M^n)$ is a solution of \eqref{eq:INT Shrinker} for some $\alpha \in [\frac{1}{n+2},\frac{1}{2}]$. We show that  $Z$ satisfies an inequality of the form 
\[\alpha K^{\alpha} b^{ij}\nabla_i \nabla_j Z + (2\alpha -1)b^{ij}\nabla_i K^\alpha \nabla_j Z \geq 0.\]  
The strong maximum principle then implies that $Z$ is constant. By examining the case of equality, we are able to show that either $\nabla h = 0$ or $\alpha = \frac{1}{n+2}$ and the cubic form vanishes. This shows that either $\Sigma$ is a round sphere, or $\alpha=\frac{1}{n+2}$ and $\Sigma$ is an ellipsoid.


Finally, we consider the case $\alpha \in (\frac{1}{2},\infty)$. As in \cite{Choi}, we consider the quantity 
\[W = K^\alpha \lambda_1^{-1} - \frac{n\alpha-1}{2n\alpha} \, |F|^2.\] 
By applying the maximum principle, we can show that any point where $W$ attains its maximum is umbilic. From this, we deduce that any maximum point of $W$ is also a maximum point of $Z$. Applying the strong maximum principle to $Z$, we are able to show that $Z$ and $W$ are both constant. This implies that $\Sigma$ is a round sphere.

\section{Preliminaries}

We first recall the notation:
\begin{itemize}
\item The metric is given by $g_{ij} = \langle F_i, F_j \rangle$, where $F_i := \nabla_i F$, and  its inverse matrix  $g^{ij}$ of $g_{ij}$, namely $g^{ij}g_{jk}=\delta^i_k$. Also, we use the notation $F^i=g^{ij} \, F_j$. 

\item We denote by $H$ and $h_{ij}$ the mean curvature and second fundamental form, respectively.

\item For a strictly convex hypersurface, we denote by $b^{ij}$ the inverse of the second fundamental form $h_{ij}$, so that $b^{ij}h_{jk}=\delta^i_k$. Moreover, $\text{\rm tr}(b)$ will denote the trace of $b$, i.e. the reciprocal of the harmonic mean curvature.

\item We denote by $\mathcal{L}$ the operator $\mathcal{L} = \alpha K^{\alpha} b^{ij}\nabla_i \nabla_j$. 

\item We denote by $C_{ijk}$ the cubic form
\begin{align*}
C_{ijk} 
&= \frac{1}{2} \, K^{-\frac{1}{n+2}} \, \nabla_k h_{ij} \\ 
&+ \frac{1}{2}  h_{jk} \nabla_i K^{-\frac{1}{n+2}}  + \frac{1}{2} h_{ki}\nabla_jK^{-\frac{1}{n+2}} + \frac{1}{2} h_{ij}\nabla_k K^{-\frac{1}{n+2}}.
\end{align*}
\end{itemize}

We next derive some basic equations.

\begin{proposition}
Given a strictly convex smooth solution $F:M^n \to \mathbb{R}^{n+1}$ of \eqref{eq:INT Shrinker}, the following equations hold:
\begin{align}
\nabla_i b^{jk} &= -b^{jl}b^{km}\nabla_i h_{lm} , \label{eq:Pre Db} \\
\mathcal{L} |F|^2 &= 2\alpha K^{\alpha} b^{ij}(g_{ij}-h_{ij}K^{\alpha})=2\alpha K^{\alpha} \, \text{\rm tr}(b) -2n\alpha K^{2\alpha}, \label{eq:Pre L|F|^2} \\
\nabla_i K^{\alpha}&=h_{ij}\langle F, F^j \rangle , \label{eq:Pre DK^a} \\
\mathcal{L} K^\alpha &=\langle F,F_i \rangle \, \nabla_i K^\alpha+n\alpha K^\alpha -\alpha K^{2\alpha} H, \label{eq:Pre LK^a} \\
\mathcal{L} h_{ij} &= -K^{-\alpha} \nabla_i K^\alpha \nabla_j K^\alpha + \alpha K^\alpha b^{pr} b^{qs} \nabla_i h_{rs} \nabla_j h_{pq} \label{eq:Pre Lh} \\ 
&+ \langle F,F_k \rangle \, \nabla^k h_{ij} + h_{ij} + (n\alpha-1) h_{ik} h_j^k K^\alpha - \alpha K^\alpha H h_{ij}, \notag \\  
\mathcal{L} b^{pq}&=  K^{-\alpha}b^{pr}b^{qs}\nabla_r K^\alpha \nabla_s K^\alpha+\alpha K^\alpha b^{pr}b^{qs}b^{ij}b^{km}\nabla_r h_{ik}\nabla_s h_{jm} \label{eq:Pre Lb} \\
 &+\langle F,F_i \rangle \, \nabla^i  b^{pq} - b^{pq} -(n\alpha-1)g^{pq} K^\alpha+\alpha K^\alpha H b^{pq} \notag.
\end{align}
\end{proposition}

\textbf{Proof.} The relation $\nabla_i (b^{jk}h_{kl})=\nabla_i \delta^j_l =0$ gives $h_{kl}\nabla_i b^{jk}=-b^{jk} \nabla_ih_{kl}$. This directly implies \eqref{eq:Pre Db}. We next compute 
\[\nabla_i \nabla_j |F|^2 = 2 \langle \nabla_i F,\nabla_j F \rangle+2 \langle F,\nabla_i \nabla_j F \rangle = 2 g_{ij} - 2 h_{ij} \langle F,\nu \rangle = 2 g_{ij} - 2 K^\alpha h_{ij},\] 
hence 
\[\mathcal{L} |F|^2 = 2\alpha K^\alpha \, \text{\rm tr}(b) - 2n \alpha K^{2\alpha}.\] 
This proves \eqref{eq:Pre L|F|^2}. 

To derive the equation \eqref{eq:Pre DK^a}, we differentiate \eqref{eq:INT Shrinker}:
\begin{align*}
\nabla_i K^{\alpha}=h_{ik}\langle F, F^k \rangle.
\end{align*}
If we differentiate this equation again, we obtain 
\begin{align*}
\nabla_i\nabla_j K^{\alpha} 
&= \langle F, F^k \rangle \, \nabla_i h_{jk}+h_{ij}-h_{ik}h^k_j\langle F, \nu \rangle\\
&= \langle F,F^k \rangle \, \nabla_k h_{ij} +h_{ij}-K^\alpha h_{ik}h^k_j, 
\end{align*} 
hence 
\[\mathcal{L} K^\alpha = \langle F,F^k \rangle \, \nabla_k K^\alpha +n\alpha K^\alpha-\alpha K^{2\alpha} H.\] 
On the other hand, using \eqref{eq:Pre Db} we compute
\begin{align*}
&\nabla_i\nabla_j K^{\alpha}=\nabla_i(\alpha K^{\alpha}b^{pq}\nabla_j h_{pq})\\
&=\alpha K^\alpha b^{pq}\nabla_i\nabla_j h_{pq}+ \alpha^2K^{\alpha}b^{rs}b^{pq}\nabla_ih_{rs} \nabla_j h_{pq}-\alpha K^\alpha b^{pr}b^{qs}\nabla_i h_{rs}\nabla_j h_{pq}.
\end{align*}
Using the commutator identity 
\begin{align*}
\nabla_i \nabla_j h_{pq}&=\nabla_i \nabla_p h_{jq}=\nabla_p \nabla_i h_{jq}+R_{ipjm}h^m_q+R_{ipqm}h^m_j \\
&=\nabla_p \nabla_q h_{ij}+(h_{ij}h_{pm}-h_{im}h_{jp})h^m_q+(h_{iq}h_{pm}-h_{im}h_{pq})h^m_j
\end{align*}
we deduce that 
\begin{align*}
\alpha K^{\alpha} b^{pq}\nabla_i \nabla_j h_{pq}&=\alpha K^{\alpha}b^{pq}\nabla_p \nabla_q h_{ij}+\alpha K^{\alpha} H h_{ij}-n\alpha K^{\alpha} h_{im}h^m_j\\
&=\mathcal{L} h_{ij}+\alpha K^{\alpha} H h_{ij}-n\alpha K^{\alpha} h_{im}h^m_j.
\end{align*} 
Combining the equations above yields 
\begin{align*}
\mathcal{L} h_{ij} = &-\alpha^2K^{\alpha}b^{rs}b^{pq}\nabla_ih_{rs} \nabla_j h_{pq}+\alpha K^\alpha b^{pr}b^{qs}\nabla_i h_{rs}\nabla_j h_{pq}\\
&+\langle F,F^k \rangle \, \nabla_k h_{ij} +h_{ij}+(n\alpha-1)h_{ik}h^k_jK^\alpha -\alpha K^{\alpha} H h_{ij}.
\end{align*} 
This completes the proof of \eqref{eq:Pre Lh}.

Finally, using \eqref{eq:Pre Db}, we obtain 
\begin{align*} 
\mathcal{L} b^{pq} 
&= \alpha K^\alpha b^{ij} \nabla_i (-b^{pr}b^{qs}\nabla_j h_{rs}) \\
&=2\alpha K^\alpha b^{ij}b^{pk}b^{rm}b^{qs}\nabla_i h_{km}\nabla_j h_{rs}-b^{pr}b^{qs} \mathcal{L} h_{rs}. 
\end{align*} 
Applying \eqref{eq:Pre Lh}, we conclude that 
\begin{align*}
\mathcal{L} b^{pq} 
&= \alpha^2 K^{\alpha}b^{pr}b^{qs}b^{ij}b^{km}\nabla_r h_{ij} \nabla_s h_{km}+\alpha K^\alpha b^{pr}b^{qs}b^{ij}b^{km}\nabla_r h_{ik}\nabla_s h_{jm}\\
&+ \langle F,F^k \rangle \, \nabla_k  b^{pq} - b^{pq} -(n\alpha-1)g^{pq} K^\alpha+\alpha K^\alpha H b^{pq}. 
\end{align*} 
Since $\nabla K^\alpha=\alpha K^\alpha b^{ij}\nabla h_{ij}$, the identity \eqref{eq:Pre Lb} follows.

\section{Classification of self-similar solutions: The case $\alpha \in [\frac{1}{n+2},\frac{1}{2}]$}

In this section, we consider the case $\alpha \in [\frac{1}{n+2},\frac{1}{2}]$. We begin with an algebraic lemma:

\begin{lemma}
\label{3rd.order.algebra}
Assume $\alpha \in [\frac{1}{n+2},\frac{1}{2}]$, $\lambda_1,\hdots,\lambda_n$ are positive real numbers (not necessarily arranged in increasing order), and $\sigma_1,\hdots,\sigma_n$ are arbitrary real numbers. Then 
\begin{align*}
Q &:= \sum_{i=1}^n \sigma_i^2 + 2 \sum_{i=1}^{n-1} \lambda_n \lambda_i^{-1} \sigma_i^2\\
&- 4\alpha \lambda_n \bigg [ \sum_{i=1}^n \lambda_i^{-1}\sigma_i+\Big((n\alpha-1)\lambda_n^{-1}-\alpha\sum^n_{i=1}\lambda_i^{-1}\Big)\sum_{i=1}^n \sigma_i \bigg ] \Big(\sum_{i=1}^n \sigma_i \Big) \\
& - 2\alpha^2 \lambda_n \, \Big(\sum_{i=1}^n \lambda_i^{-1} \Big)  \Big(\sum_{i=1}^n \sigma_i \Big)^2 +\big[2n\alpha^2 + (n-1)\alpha - 1\big] \, \Big(\sum_{i=1}^n \sigma_i \Big)^2 \\ 
&\geq 0. 
\end{align*}
Moreover, if equality holds, then we either have $\sigma_1 = \hdots = \sigma_n = 0$ or we have $\alpha = \frac{1}{n+2}$ and $\sigma_1=\hdots=\sigma_{n-1}=\frac{1}{3} \, \sigma_n$.
\end{lemma}

\textbf{Proof.} 
If $\sum_{i=1}^n \sigma_i = 0$, the assertion is trivial. Hence, it suffices to consider the case $\sum_{i=1}^n \sigma_i \neq 0$. By scaling, we may assume $\sum_{i=1}^n \sigma_i = 1$. Let us define 
\[\tau_i = \begin{cases} \sigma_i - \alpha & \text{\rm for $i = 1,\hdots,n-1$} \\ \sigma_n - 1 + (n-1)\alpha & \text{\rm for $i=n$.} \end{cases}\] 
Then 
\[\sum_{i=1}^n \tau_i = \sum_{i=1}^n \sigma_i - 1 = 0\] 
and 
\[\sum_{i=1}^n \lambda_i^{-1} \tau_i = \sum_{i=1}^n \lambda_i^{-1}\sigma_i+(n\alpha-1)\lambda_n^{-1}-\alpha\sum^n_{i=1}\lambda_i^{-1}.\] 
Therefore, the quantity $Q$ satisfies
\begin{align*} 
Q &= \sum_{i=1}^{n-1} (\tau_i+ \alpha)^2+\big(\tau_n+1-(n-1)\alpha\big)^2 + 2 \sum_{i=1}^{n-1} \lambda_n \lambda_i^{-1} (\tau_i+\alpha)^2\\
&-4\alpha  \sum_{i=1}^n \lambda_n \lambda_i^{-1}\tau_i  - 2\alpha^2 \sum_{i=1}^n \lambda_n \lambda_i^{-1}  + 2n\alpha^2 + (n-1)\alpha - 1\\
&= \sum_{i=1}^n \tau_i^2+ 2\alpha \sum_{i=1}^{n-1} \tau_i + 2 (1-(n-1)\alpha) \tau_n + (n-1)\alpha^2 + (1-(n-1)\alpha)^2\\
& + 2 \sum_{i=1}^{n-1} \lambda_n \lambda_i^{-1} (\tau_i^2+2\alpha \tau_i+\alpha^2)-4\alpha \sum_{i=1}^{n-1} \lambda_n \lambda_i^{-1} \tau_i-4\alpha \tau_n\\ 
&- 2\alpha^2  \sum_{i=1}^{n-1} \lambda_n \lambda_i^{-1} - 2\alpha^2 + 2n\alpha^2 + (n-1)\alpha - 1 \\ 
&= \sum_{i=1}^n \tau_i^2+ 2\alpha \sum^{n-1}_{i=1}\tau_i  + 2 (1-(n+1)\alpha) \tau_n + 2 \sum_{i=1}^{n-1} \lambda_n \lambda_i^{-1}\tau_i^2 \\ 
&+ (n-1)\alpha((n+2)\alpha-1). 
\end{align*} 
Using the identity $\sum_{i=1}^n \tau_i = 0$, we obtain 
\begin{align*}
Q = \sum_{i=1}^n \tau_i^2 + 2 (1-(n+2)\alpha) \tau_n + 2 \sum_{i=1}^{n-1} \lambda_n \lambda_i^{-1}\tau_i^2+ (n-1)\alpha ((n+2)\alpha-1).
\end{align*}
Moreover, the identity $\sum^n_{i=1}\tau_i=0$ gives
\begin{align*}
\sum_{i=1}^n \tau_i^2 &= \sum_{i=1}^{n-1} \Big(\tau_i+\frac{1}{n-1}\tau_n \Big)^2-\frac{2}{n-1}\tau_n \sum^{n-1}_{i=1}\tau_i+\frac{n-2}{n-1}\tau_n^2\\
&= \sum_{i=1}^{n-1} \Big(\tau_i+\frac{1}{n-1}\tau_n \Big)^2+\frac{n}{n-1}\tau_n^2.
\end{align*}
Thus, 
\begin{align*}
Q &= \frac{n}{n-1} \, \Big (\tau_n + \frac{n-1}{n} \, (1-(n+2)\alpha) \Big )^2 + \sum_{i=1}^{n-1} \Big(\tau_i + \frac{1}{n-1} \, \tau_n\Big)^2 \\ 
& + 2 \sum_{i=1}^{n-1} \lambda_n \lambda_i^{-1} \tau_i^2+ \frac{n-1}{n} \, (1-2\alpha) ((n+2)\alpha-1). 
\end{align*}
The right hand side is clearly nonnegative. Moreover, if equality holds, then $\tau_1=\hdots=\tau_{n}=0$, and $\alpha = \frac{1}{n+2}$. This proves the lemma. \\

\begin{theorem}\label{soliton.uniqueness}
Assume $\alpha \in [\frac{1}{n+2},\frac{1}{2}]$ and $\Sigma$ is a strictly convex closed smooth solution of \eqref{eq:INT Shrinker}. Then, either $\Sigma$ is a round sphere, or $\alpha=\frac{1}{n+2}$ and $\Sigma$ is an ellipsoid.
\end{theorem}

\textbf{Proof.} 
Taking the trace in equation \eqref{eq:Pre Lb} gives 
\begin{align*}
\mathcal{L} \, \text{\rm tr}(b) 
&=  K^{-\alpha}b^{pr}b_p^s\nabla_r K^\alpha \nabla_s K^\alpha+\alpha K^\alpha b^{pr}b_p^sb^{ij}b^{km}\nabla_r h_{ik}\nabla_s h_{jm} \\ 
&+\langle F,F_i \rangle \, \nabla^i  \text{\rm tr}(b) - \text{\rm tr}(b) - n(n\alpha-1) K^\alpha+\alpha K^\alpha H \text{\rm tr}(b).
\end{align*} 
Using equation \eqref{eq:Pre LK^a}, we obtain 
\begin{align*} 
\mathcal{L} (K^\alpha \, \text{\rm tr}(b)) 
&= K^\alpha \mathcal{L} \, \text{\rm tr}(b) + \mathcal{L} K^\alpha \, \text{\rm tr}(b) + 2\alpha K^\alpha b^{ij} \nabla_i K^\alpha \nabla_j \text{\rm tr}(b) \\ 
&=  b^{pr}b_p^s\nabla_r K^\alpha \nabla_s K^\alpha+\alpha K^{2\alpha} b^{pr}b_p^sb^{ij}b^{km}\nabla_r h_{ik}\nabla_s h_{jm} \\ 
&+ \langle F,F_i \rangle \, \nabla^i(K^\alpha \, \text{\rm tr}(b)) + (n\alpha-1) K^\alpha \, \text{\rm tr}(b) - n(n\alpha-1) K^{2\alpha} \\ 
&+ 2\alpha K^\alpha b^{ij} \nabla_i K^\alpha \nabla_j \text{\rm tr}(b). 
\end{align*} 
Using \eqref{eq:Pre L|F|^2}, it follows that the function 
\begin{equation}
Z=K^\alpha\, \text{\rm tr}(b)-\frac{n\alpha-1}{2\alpha} |F|^2 \label{eq:Z}
\end{equation} 
satisfies 
\begin{align}
\mathcal{L} Z &=  b^{pr}b_p^s\nabla_r K^\alpha \nabla_s K^\alpha+\alpha K^{2\alpha} b^{pr}b_p^sb^{ij}b^{km}\nabla_r h_{ik}\nabla_s h_{jm} \notag \\ 
&+ \langle F,F_i \rangle \, \nabla^i(K^\alpha \, \text{\rm tr}(b)) + 2\alpha K^\alpha b^{ij} \nabla_i K^\alpha \nabla_j \text{\rm tr}(b). \label{eq:LZ}
\end{align}
Using \eqref{eq:Pre DK^a}, we obtain 
\[\frac{1}{2}\nabla^i |F|^2=\langle F,F^i\rangle=b^{ij}\nabla_j K^\alpha,\] 
hence 
\begin{align*}
\nabla_i Z=K^\alpha\nabla_i \text{\rm tr}(b) +\text{\rm tr}(b) \nabla_i K^\alpha-\frac{n\alpha-1}{\alpha} b_i^j\nabla_j K^\alpha.
\end{align*}
This gives 
\begin{align*} 
\langle F,F_i \rangle \, \nabla^i(K^\alpha \, \text{\rm tr}(b)) 
&= b^{ij} \nabla_i K^\alpha \nabla_j Z + \frac{n\alpha-1}{\alpha} b^{ik} b_k^j \nabla_i K^\alpha \nabla_j K^\alpha 
\end{align*}
and 
\begin{align*} 
2\alpha K^\alpha b^{ij} \nabla_i K^\alpha \nabla_j \text{\rm tr}(b) 
&= 2\alpha b^{ij}\nabla_i K^\alpha \nabla_j Z \\ 
&-\big(2\alpha b^{ij}\text{\rm tr}(b) -2(n\alpha-1)b^{ik}b_k^j\big) \nabla_i K^\alpha \nabla_j K^\alpha.
\end{align*}
Substituting these identities into \eqref{eq:LZ} gives 
\begin{align}\label{eq:Z} 
&\mathcal{L} Z + (2\alpha-1) b^{ij} \nabla_i K^\alpha \nabla_j Z \\ 
&= 4\alpha b^{ij}\nabla_i K^\alpha \nabla_j Z +\alpha K^{2\alpha} b^{pr}b_p^sb^{ij}b^{km}\nabla_r h_{ik}\nabla_s h_{jm} \notag \\
&+\big(-2\alpha b^{ij}\text{\rm tr}(b) +(2n\alpha+n-1-\alpha^{-1})b^{ik}b_k^j\big) \nabla_i K^\alpha \nabla_j K^\alpha. \notag
\end{align}
Let us fix an arbitrary point $p$. We can choose an orthonormal frame so that $h_{ij}(p)=\lambda_i \delta_{ij}$. With this understood, we have 
\begin{align}\label{eq:Max DK^a}
&\nabla_i K^\alpha=\alpha K^\alpha \sum_{j=1}^n\lambda_j^{-1}\nabla_i h_{jj}, &&\nabla_i \text{\rm tr}(b)=-\sum^n_{j=1}\lambda_j^{-2}\nabla_i h_{jj}.
\end{align}
Let $D$ denote the set of all triplets $(i,j,k)$ such that $i,j,k$ are pairwise distinct. Then, by using \eqref{eq:Z} and \eqref{eq:Max DK^a}, we have
\begin{align}\label{eq:Max Final eq}
&\alpha^{-1}K^{-2\alpha}\Big(\mathcal{L} Z+ (2\alpha -1)b^{ij}\nabla_i K^\alpha \nabla_j Z\Big)\\  
&=\sum_D \lambda_i^{-2}\lambda_j^{-1}\lambda_k^{-1} (\nabla_i h_{jk})^2 +4\alpha  \sum_k\lambda^{-1}_k (\nabla_k \log K) (K^{-\alpha}\nabla_k Z) \notag \\
&+\sum_k \sum_i \lambda_k^{-2} \lambda_i^{-2} \, (\nabla_k h_{ii})^2 + 2 \sum_k \sum_{i \neq k} \lambda_k^{-1} \lambda_i^{-3} \, (\nabla_k h_{ii})^2 \notag \\ 
&+\sum_k \lambda_k^{-1} \, \Big [ -2\alpha^2 \, \text{\rm tr}(b)  + \big(2n\alpha^2 + (n-1)\alpha - 1\big) \lambda_k^{-1} \Big ] \, (\nabla_k \log K)^2. \notag
\end{align} 
We claim that for each $k$ the following holds:
\begin{align}\label{eq:gradient terms}
&\sum_i \lambda_k^{-1} \lambda_i^{-2} \, (\nabla_k h_{ii})^2 + 2 \sum_{i \neq k}  \lambda_i^{-3} \, (\nabla_k h_{ii})^2 +4\alpha (\nabla_k \log K) (K^{-\alpha}\nabla_k Z)\\ 
&+ \Bigl[- 2\alpha^2 \, \text{\rm tr}(b)  + \big(2n\alpha^2 + (n-1)\alpha - 1\big) \lambda_k^{-1} \Bigr] \, (\nabla_k \log K)^2 \geq 0. \notag
\end{align} 
Notice that
\[K^{-\alpha} \, \nabla_k Z = -\sum_i \lambda_i^{-2}\nabla_k h_{ii} + (\alpha \, \text{\rm tr}(b) - (n\alpha-1) \, \lambda_k^{-1}) \, \nabla_k \log K.\] 
After relabeling indices, we may assume $k=n$. If we put $\sigma_i := \lambda_i^{-1} \, \nabla_n h_{ii}$, then the assertion follows from Lemma \ref{3rd.order.algebra}. This proves \eqref{eq:gradient terms}. Combining \eqref{eq:Max Final eq} and \eqref{eq:gradient terms}, we conclude that 
\[\mathcal{L} Z+ (2\alpha -1)b^{ij}\nabla_i K^\alpha \nabla_j Z \geq 0\] 
at each point $p$. Therefore, by the strong maximum principle, $Z$ is a constant. Hence, the left hand side of \eqref{eq:Max Final eq} is zero. Therefore, $\nabla_i h_{jk} = 0$ if $i,j,k$ are all distinct. Moreover, since we have equality in the lemma, we either have $\lambda_i^{-1} \, \nabla_k h_{ii} = 0$ for all $i,k$, or we have $\alpha=\frac{1}{n+2}$ and $\lambda_i^{-1} \, \nabla_k h_{ii} = \frac{1}{3} \, \lambda_k^{-1} \, \nabla_k h_{kk}$ for $i \neq k$. 

In the first case, we conclude that $\nabla_i h_{jk} = 0$ for all $i,j,k$ and thus $\Sigma$ is a round sphere. 

In the second case, we obtain $\lambda_i^{-1} \, \nabla_k h_{ii} = \frac{1}{n+2} \, \nabla_k \log K$ for $i \neq k$ and $\lambda_k^{-1} \, \nabla_k h_{kk} = \frac{3}{n+2} \, \nabla_k \log K$. This gives 
\begin{align*} 
C_{ijk} 
&= \frac{1}{2} \, K^{-\frac{1}{n+2}} \, \nabla_k h_{ij} \\ 
&+ \frac{1}{2} \, h_{jk}\nabla_iK^{-\frac{1}{n+2}}  + \frac{1}{2} h_{ki} \nabla_jK^{-\frac{1}{n+2}} + \frac{1}{2} h_{ij} \nabla_kK^{-\frac{1}{n+2}} \\ 
&= 0
\end{align*}
for all $i,j,k$. Since the cubic form $C_{ijk}$ vanishes everywhere, the surface is an ellipsoid by the Berwald-Pick theorem (see e.g. \cite{Nomizu-Sasaki}, Theorem 4.5). This proves the theorem.\\

\section{Classification of self-similar solutions: The case $\alpha \in (\frac{1}{2},\infty)$}

We now turn to the case $\alpha \in (\frac{1}{2},\infty)$. In the following, we denote by $\lambda_1 \leq \hdots \leq \lambda_n$ the eigenvalues of the second fundamental form, arranged in increasing order. Each eigenvalue defines a Lipschitz continuous function on $M$.

\begin{lemma}
Suppose that $\varphi$ is a smooth function such that $\lambda_1 \geq \varphi$ everywhere and $\lambda_1 = \varphi$ at $\bar{p}$. Let $\mu$ denote the multiplicity of the smallest curvature eigenvalue at $\bar{p}$, so that $\lambda_1 = \hdots = \lambda_\mu < \lambda_{\mu+1} \leq \hdots \leq \lambda_n$. Then, at $\bar{p}$, $\nabla_i h_{kl} = \nabla_i \varphi \, \delta_{kl}$ for $1 \leq k,l \leq \mu$. Moreover, 
\[\nabla_i \nabla_i \varphi \leq \nabla_i \nabla_i h_{11}-2 \sum_{l>\mu} (\lambda_l-\lambda_1)^{-1} \, (\nabla_i h_{1l})^2.\] 
at $\bar{p}$.
\end{lemma}

\textbf{Proof.} 
Fix an index $i$, and let $\gamma(s)$ be the geodesic satisfying $\gamma(0)=\bar{p}$ and $\gamma'(0)=e_i$. Moreover, let $v(s)$ be a vector field along $\gamma$ such that $v(0) \in \text{\rm span}\{e_1,\hdots,e_\mu\}$ and $v'(0) \in \text{\rm span}\{e_{\mu+1},\hdots,e_n\}$. Then the function $s \mapsto h(v(s),v(s)) - \varphi(\gamma(s)) \, |v(s)|^2$ has a local minimum at $s=0$. This gives 
\begin{align*} 
0 
&= \frac{d}{ds} \big ( h(v(s),v(s)) - \varphi(\gamma(s)) \, |v(s)|^2 \big ) \Big |_{s=0} \\ 
&= \nabla_i h(v(0),v(0)) + 2 \, h(v(0),v'(0)) - \nabla_i \varphi \, |v(0)|^2 - 2 \, \langle v(0),v'(0) \rangle \\ 
&=  \nabla_i h(v(0),v(0)) - \nabla_i \varphi \, |v(0)|^2. 
\end{align*}
Since $v(0) \in \text{\rm span}\{e_1,\hdots,e_\mu\}$ is arbitrary, it follows that $\nabla_i h_{kl} = \nabla_i \varphi \, \delta_{kl}$ for $1 \leq k,l \leq \mu$ at the point $\bar{p}$.

We next consider the second derivative. To that end, we choose $v(0)=e_1$, $v'(0) = -\sum_{l>\mu} (\lambda_l-\lambda_1)^{-1} \, \nabla_i h_{1l} \, e_l$, and $v''(0)=0$. Since the function $s \mapsto h(v(s),v(s)) - \varphi(\gamma(s)) \, |v(s)|^2$ has a local minimum at $s=0$, we obtain  
\begin{align*} 
0 &\leq \frac{d^2}{ds^2} \big ( h(v(s),v(s)) - \varphi(\gamma(s)) \, |v(s)|^2 \big ) \Big |_{s=0} \\ 
&= \nabla_i \nabla_i h(v(0),v(0)) + 4 \nabla_i h(v(0),v'(0)) + 2 h(v'(0),v'(0)) \\ 
&- \nabla_i \nabla_i \varphi \, |v(0)|^2 - 4 \nabla_i \varphi \, \langle v(0),v'(0) \rangle - 2 \varphi \, |v'(0)|^2 \\ 
&= \nabla_i \nabla_i h_{11} - 4 \sum_{l>\mu} (\lambda_l-\lambda_1)^{-1} \, (\nabla_i h_{1l})^2 + 2 \sum_{l>\mu} \lambda_l (\lambda_l-\lambda_1)^{-2} \, (\nabla_i h_{1l})^2 \\ 
&- \nabla_i \nabla_i \varphi - 2 \sum_{l>\mu} \lambda_1 (\lambda_l-\lambda_1)^{-2} \, (\nabla_i h_{1l})^2 \\ 
&= \nabla_i \nabla_i h_{11} - 2 \sum_{l>\mu} (\lambda_l-\lambda_1)^{-2} \, (\nabla_i h_{1l})^2 - \nabla_i \nabla_i \varphi. 
\end{align*} 
This proves the assertion.

\begin{theorem}\label{soliton.uniqueness2}
Assume $\alpha \in (\frac{1}{2},\infty)$ and $\Sigma$ is a strictly convex closed smooth solution of \eqref{eq:INT Shrinker}. Then $\Sigma$ is a round sphere.
\end{theorem}

\textbf{Proof.} 
Let us consider the function 
\begin{equation}
W = K^\alpha \lambda_1^{-1} - \frac{n\alpha-1}{2n\alpha} \, |F|^2. \label{eq:W}
\end{equation}
Let us consider an arbitrary point $\bar{p}$ where $W$ attains its maximum. As above, we denote by $\mu$ the multiplicity of the smallest eigenvalue of the second fundamental form. Let us define a smooth function $\varphi$ such that 
\[W(\bar{p}) = K^\alpha \varphi^{-1} - \frac{n\alpha-1}{2n\alpha} \, |F|^2.\] 
Since $W$ attains its maximum at $\bar{p}$, we have $\lambda_1 \geq \varphi$ everywhere and $\lambda_1 = \varphi$ at $\bar{p}$. Therefore, we may apply the previous lemma. Hence, at the point $\bar{p}$, we have $\nabla_i h_{kl} = \nabla_i \varphi \, \delta_{kl}$ for $1 \leq k,l \leq \mu$. Moreover, 
\[\nabla_k \nabla_k \varphi \leq \nabla_k \nabla_k h_{11}-2 \sum_{l>\mu} (\lambda_l-\lambda_1)^{-1} \, (\nabla_k h_{1l})^2\] 
at $\bar{p}$. We multiply both sides by $\alpha K^\alpha \lambda_k^{-1}$ and sum over $k$. This gives 
\[\mathcal{L} \varphi \leq \mathcal{L} h_{11} - 2\alpha K^\alpha \sum_k \sum_{l>\mu} \lambda_k^{-1} (\lambda_l-\lambda_1)^{-1} \, (\nabla_k h_{1l})^2.\] 
By \eqref{eq:Pre Lh}, we have  
\begin{align*} 
\mathcal{L} h_{ij} &= -K^{-\alpha} \nabla_i K^\alpha \nabla_j K^\alpha + \alpha K^\alpha b^{pr} b^{qs} \nabla_i h_{rs} \nabla_j h_{pq} \\ 
&+ \langle F,F_k \rangle \, \nabla^k h_{ij} + h_{ij} + (n\alpha-1) h_{ik} h_j^k K^\alpha - \alpha K^\alpha H h_{ij}. 
\end{align*} 
Thus, 
\begin{align*} 
\mathcal{L} \varphi &\leq -2\alpha K^\alpha \sum_k \sum_{l>\mu} \lambda_k^{-1} (\lambda_l-\lambda_1)^{-1} \, (\nabla_k h_{1l})^2 \\ 
&-K^{-\alpha} (\nabla_1 K^\alpha)^2 + \alpha K^\alpha \sum_{k,l} \lambda_k^{-1} \lambda_l^{-1} (\nabla_k h_{1l})^2 \\ 
&+ \sum_k \langle F,F_k \rangle \, \nabla_k h_{11} + \lambda_1 + (n\alpha-1) \lambda_1^2 K^\alpha - \alpha K^\alpha H \lambda_1. 
\end{align*} 
Using the estimate 
\begin{align*} 
&-2\alpha K^\alpha \sum_k \sum_{l>\mu} \lambda_k^{-1} (\lambda_l-\lambda_1)^{-1} \, (\nabla_k h_{1l})^2 + \alpha K^\alpha \sum_{k,l} \lambda_k^{-1} \lambda_l^{-1} (\nabla_k h_{1l})^2 \\ 
&= -\alpha K^\alpha \sum_k \sum_{l>\mu} \lambda_k^{-1} [2(\lambda_l-\lambda_1)^{-1}-\lambda_l^{-1}] \, (\nabla_k h_{1l})^2 + \alpha K^\alpha \sum_k \lambda_k^{-1} \lambda_1^{-1} (\nabla_k h_{11})^2 \\ 
&\leq -\alpha K^\alpha \sum_{l>\mu} \lambda_1^{-1} [2(\lambda_l-\lambda_1)^{-1}-\lambda_l^{-1}] \, (\nabla_1 h_{1l})^2 + \alpha K^\alpha \sum_k \lambda_k^{-1} \lambda_1^{-1} (\nabla_k h_{11})^2 \\ 
&= -2\alpha K^\alpha \sum_{k>\mu} \lambda_1^{-1} [(\lambda_k-\lambda_1)^{-1}-\lambda_k^{-1}] \, (\nabla_k h_{11})^2 + \alpha K^\alpha \lambda_1^{-2} (\nabla_1 h_{11})^2, 
\end{align*}
we obtain
\begin{align*} 
\mathcal{L} \varphi &\leq -2\alpha K^\alpha \lambda_1^{-1} \sum_{k>\mu} [(\lambda_k-\lambda_1)^{-1}-\lambda_k^{-1}] \, (\nabla_k h_{11})^2 \\ 
&+ \alpha K^\alpha \lambda_1^{-2} (\nabla_1 h_{11})^2 - K^{-\alpha} (\nabla_1 K^\alpha)^2  \\ 
&+ \sum_k \langle F,F_k \rangle \, \nabla_k h_{11} + \lambda_1 + (n\alpha-1) \lambda_1^2 K^\alpha - \alpha K^\alpha H \lambda_1. 
\end{align*} 
Since $\nabla_k \varphi = \nabla_k h_{11}$, it follows that 
\begin{align*} 
\mathcal{L} (\varphi^{-1}) &\geq 2\alpha K^\alpha \lambda_1^{-3} \sum_k \lambda_k^{-1} (\nabla_k h_{11})^2 \\ 
&+ 2\alpha K^\alpha \lambda_1^{-3} \sum_{k>\mu} [(\lambda_k-\lambda_1)^{-1}-\lambda_k^{-1}] \, (\nabla_k h_{11})^2 \\ 
&- \alpha K^\alpha \lambda_1^{-4} (\nabla_1 h_{11})^2 + K^{-\alpha} \lambda_1^{-2} (\nabla_1 K^\alpha)^2 \\
&+ \sum_k \langle F,F_k \rangle \, \nabla_k (\lambda_1^{-1}) - \lambda_1^{-1} - (n\alpha-1) K^\alpha + \alpha K^\alpha H \lambda_1^{-1} \\  
&= 2\alpha K^\alpha \lambda_1^{-3} \sum_{k>\mu} (\lambda_k-\lambda_1)^{-1} \, (\nabla_k h_{11})^2 \\ 
&+ \alpha K^\alpha \lambda_1^{-4} (\nabla_1 h_{11})^2 + K^{-\alpha} \lambda_1^{-2} (\nabla_1 K^\alpha)^2 \\
&+ \sum_k \langle F,F_k \rangle \, \nabla_k (\lambda_1^{-1}) - \lambda_1^{-1} - (n\alpha-1) K^\alpha + \alpha K^\alpha H \lambda_1^{-1}. 
\end{align*} 
This gives 
\begin{align*} 
\mathcal{L} (K^\alpha \varphi^{-1}) 
&= K^\alpha \mathcal{L}(\varphi^{-1}) + \varphi^{-1} \mathcal{L}(K^\alpha) + 2\alpha K^\alpha \sum_k \lambda_k^{-1} \, \nabla_k K^\alpha \, \nabla_k (\varphi^{-1}) \\ 
&\geq 2\alpha \sum_k \lambda_k^{-1} \, \nabla_k K^\alpha \, \nabla_k (K^\alpha \varphi^{-1}) - 2\alpha \lambda_1^{-1} \sum_k \lambda_k^{-1} \, (\nabla_k K^\alpha)^2 \\ 
&+ 2\alpha K^{2\alpha} \lambda_1^{-3} \sum_{k>\mu} (\lambda_k-\lambda_1)^{-1} \, (\nabla_k h_{11})^2 \\ 
&+ \alpha K^{2\alpha} \lambda_1^{-4} (\nabla_1 h_{11})^2 + \lambda_1^{-2} (\nabla_1 K^\alpha)^2 \\
&+ \sum_k \langle F,F_k \rangle \, \nabla_k (K^\alpha \varphi^{-1}) + (n\alpha-1) K^\alpha \lambda_1^{-1} - (n\alpha-1) K^{2\alpha}. 
\end{align*} 
By assumption, the function $K^\alpha \varphi^{-1} - \frac{n\alpha-1}{2n\alpha} \, |F|^2$ is constant. Consequently, $\nabla_k (K^\alpha \varphi^{-1}) = \frac{n\alpha-1}{2n\alpha} \, \nabla_k |F|^2$, and
\begin{align*}
0 &= \mathcal{L}(K^\alpha \varphi^{-1} - \frac{n\alpha-1}{2n\alpha} \, |F|^2) \\ 
&\geq \frac{n\alpha-1}{n} \sum_k \lambda_k^{-1} \, \nabla_k K^\alpha \, \nabla_k |F|^2 - 2\alpha \lambda_1^{-1} \sum_k \lambda_k^{-1} \, (\nabla_k K^\alpha)^2 \\ 
&+ 2\alpha K^{2\alpha} \lambda_1^{-3} \sum_{k>\mu} (\lambda_k-\lambda_1)^{-1} \, (\nabla_k h_{11})^2 \\ 
&+ \alpha K^{2\alpha} \lambda_1^{-4} (\nabla_1 h_{11})^2 + \lambda_1^{-2} (\nabla_1 K^\alpha)^2 \\
&+ \frac{n\alpha-1}{2n\alpha} \sum_k \langle F,F_k \rangle \, \nabla_k |F|^2 + (n\alpha-1) K^\alpha (\lambda_1^{-1} - \frac{1}{n} \, \text{\rm tr}(b)). 
\end{align*} 
Recall that 
\[\frac{1}{2} \, \nabla_k |F|^2 = \langle F,F_k \rangle = \lambda_k^{-1} \nabla_k K^\alpha\] 
by \eqref{eq:Pre DK^a}. Moreover, using the identity $\nabla_k \varphi = \nabla_k h_{11}$, we obtain   
\begin{align*} 
0 &= \nabla_k (K^\alpha \varphi^{-1}) - \frac{n\alpha-1}{2n\alpha} \, \nabla_k |F|^2 \\ 
&= (\lambda_1^{-1} - \frac{n\alpha-1}{n\alpha} \lambda_k^{-1}) \nabla_k K^\alpha - K^\alpha \lambda_1^{-2} \nabla_k h_{11}
\end{align*} 
at $\bar{p}$. Note that, if $2 \leq l \leq \mu$, then $\nabla_k h_{1l} = 0$ for all $k$. Putting $k=1$ gives $\nabla_l h_{11}=0$ and $\nabla_l K = 0$ for $2 \leq l \leq \mu$. Putting these facts together, we obtain 
\begin{align*} 
0 &\geq \frac{2(n\alpha-1)}{n} \sum_k \lambda_k^{-2} \, (\nabla_k K^\alpha)^2 - 2\alpha \lambda_1^{-1} \sum_k \lambda_k^{-1} \, (\nabla_k K^\alpha)^2 \\ 
&+ 2\alpha \lambda_1 \sum_{k>\mu} (\lambda_k-\lambda_1)^{-1} \, (\lambda_1^{-1} - \frac{n\alpha-1}{n\alpha} \lambda_k^{-1})^2 \, (\nabla_k K^\alpha)^2 \\ 
&+ \frac{1}{n^2\alpha} \, \lambda_1^{-2} \, (\nabla_1 K^\alpha)^2 + \lambda_1^{-2} (\nabla_1 K^\alpha)^2 \\
&+ \frac{n\alpha-1}{n\alpha} \sum_k \lambda_k^{-2} \, (\nabla_k K^\alpha)^2 + (n\alpha-1) K^\alpha (\lambda_1^{-1} - \frac{1}{n} \, \text{\rm tr}(b)). 
\end{align*} 
Using the identities 
\begin{align*} 
&\frac{2(n\alpha-1)}{n} \, \lambda_k^{-2} - 2\alpha \lambda_1^{-1} \lambda_k^{-1} \\ 
&+ 2\alpha \lambda_1 (\lambda_k-\lambda_1)^{-1} \, (\lambda_1^{-1} - \frac{n\alpha-1}{n\alpha} \lambda_k^{-1})^2 \\ 
&+ \frac{n\alpha-1}{n\alpha} \, \lambda_k^{-2} \\ 
&= \Big ( \frac{n\alpha-1}{n\alpha} + \frac{2}{n} + \frac{2}{n^2\alpha} \, \lambda_1 (\lambda_k-\lambda_1)^{-1} \Big ) \, \lambda_k^{-2}
\end{align*} 
and 
\begin{align*} 
&\frac{2(n\alpha-1)}{n} \, \lambda_1^{-2} - 2\alpha \lambda_1^{-2} + \frac{1}{n^2\alpha} \, \lambda_1^{-2} + \lambda_1^{-2} + \frac{n\alpha-1}{n\alpha} \, \lambda_1^{-2} \\ 
&= \frac{n-1}{n^2\alpha} \, (2n\alpha-1) \, \lambda_1^{-2}, 
\end{align*}
the previous inequality can be rewritten as follows: 
\begin{align*} 
0 &\geq \sum_{k>\mu} \Big ( \frac{n\alpha-1}{n\alpha} + \frac{2}{n} + \frac{2}{n^2\alpha} \, \lambda_1 (\lambda_k-\lambda_1)^{-1} \Big ) \, \lambda_k^{-2} \, (\nabla_k K^\alpha)^2 \\ 
&+ \frac{n-1}{n^2\alpha} \, (2n\alpha-1) \, \lambda_1^{-2} \, (\nabla_1 K^\alpha)^2 + (n\alpha-1) K^\alpha (\lambda_1^{-1} - \frac{1}{n} \, \text{\rm tr}(b)). 
\end{align*} 
Since $\alpha > \frac{1}{n}$, it follows that $\bar{p}$ is an umbilic point. Since $\bar{p}$ is an umbilic point and $\alpha > \frac{1}{2}$, there exists a neighborhood $U$ of $\bar{p}$ with the property that 
\begin{align*} 
&\alpha^{-1}K^{-2\alpha}\Big(\mathcal{L} Z - (2\alpha +1)b^{ij}\nabla_i K^\alpha \nabla_j Z\Big)\\ 
&= \sum_D \lambda_i^{-2}\lambda_j^{-1}\lambda_k^{-1} (\nabla_i h_{jk})^2 \\
&+\sum_k \sum_i \lambda_k^{-2} \lambda_i^{-2} \, (\nabla_k h_{ii})^2 + 2 \sum_k \sum_{i \neq k} \lambda_k^{-1} \lambda_i^{-3} \, (\nabla_k h_{ii})^2 \\ 
&+\sum_k \lambda_k^{-1} \, \Big [ -2\alpha^2 \, \text{\rm tr}(b)  + \big(2n\alpha^2 + (n-1)\alpha - 1\big) \lambda_k^{-1} \Big ] \, (\nabla_k \log K)^2 \\ 
&\geq 0
\end{align*} 
at each point in $U$. (Indeed, if $n \geq 3$, the last inequality follows immediately from the fact that $(n-1)\alpha-1 \geq 0$. For $n=2$ the last inequality follows from a straightforward calculation.)

Now, since $\bar{p}$ is an umbilic point, we have $Z(p) \leq nW(p) \leq nW(\bar{p}) = Z(\bar{p})$ for each point $p \in U$. Thus, $Z$ attains a local maximum at $\bar{p}$. By the strong maximum principle, $Z(p) = Z(\bar{p})$ for all points $p \in U$. This implies $W(p)=W(\bar{p})$ for all points $p \in U$. Thus, the set of all points where $W$ attains its maximum is open. Consequently, $W$ is constant. This implies that $\Sigma$ is umbilic, hence a round sphere.

\section{Proof of Theorem \ref{main.thm}}

Suppose we have any strictly convex solution to the flow with speed $-K^\alpha \nu$, where $\alpha \in [\frac{1}{n+2},\infty)$. It is known that the flow converges to a soliton after rescaling; for $\alpha>\frac{1}{n+2}$, this follows from Theorem 6.2 in \cite{Andrews-Guan-Ni}, while for $\alpha=\frac{1}{n+2}$ this follows from results in Section 9 of \cite{Andrews3}. By Theorem \ref{soliton.uniqueness} and Theorem \ref{soliton.uniqueness2}, either the limit is a round sphere, or $\alpha=\frac{1}{n+2}$ and the limit is an ellipsoid.

\end{document}